\documentclass[12pt,oneside,reqno]{amsart}

\usepackage{amsmath}
\usepackage{amssymb}
\usepackage{graphics}

\renewcommand{\bar}{\overline}

\newcommand{\eps}{\varepsilon}

\newcommand{\CC}{\mathbb{C}}

\newcommand{\FF}{\mathbb{F}}

\newcommand{\PP}{\mathbb{P}}
\newcommand{\QQ}{\mathbb{Q}}

\newcommand{\ZZ}{\mathbb{Z}}

\newcommand{\Qp}{\QQ_p}

\newcommand{\Cp}{\CC_p}

\newcommand{\CK}{\CC_K}
\newcommand{\Qpbar}{\bar{\QQ}_p}

\newcommand{\Fpbar}{\bar{\FF}_p}

\newcommand{\ints}{{\mathcal O}}

\newcommand{\maxid}{{\mathcal M}}

\newcommand{\khat}{\hat{k}}

\newcommand{\Fat}{{\mathcal F}}
\newcommand{\Jul}{{\mathcal J}}

\newcommand{\PCK}{\PP^1(\CK)}

\newcommand{\PC}{\PP^1(\CC)}

\newcommand{\PK}{\PP^1(K)}

\newcommand{\rad}{\hbox{\rm rad}}

\DeclareMathOperator{\charact}{char}

\newcommand{\Dbar}{\bar{D}}

\newcommand{\dsps}{\displaystyle}

\theoremstyle{plain}
\newtheorem{thm}{Theorem}[section]

\newtheorem{lemma}[thm]{Lemma}

\newtheorem*{thmblank}{Theorem}

\theoremstyle{remark}
\newtheorem{remark}[thm]{Remark}

\title{Wandering Domains in Non-Archimedean Polynomial Dynamics}
\author{Robert L. Benedetto}
\date{November 8, 2003; revised February 10, 2006}
\address{Department of Mathematics and Computer Science \\
        Amherst College \\
        Amherst, MA 01002 \\
        USA}
\email{rlb@cs.amherst.edu}
\urladdr{http://www.cs.amherst.edu/\textasciitilde rlb}
\thanks{The author gratefully acknowledges the support
        of NSF grant DMS-0071541.}

\begin{document}

\begin{abstract}
We extend a recent
result on the existence of wandering domains
of polynomial functions defined over the $p$-adic
field $\Cp$ to
any algebraically closed complete non-archimedean field $\CK$
with residue characteristic $p>0$.
We also prove that polynomials with wandering domains form
a dense subset of a certain one-dimensional family
of degree $p+1$ polynomials in $\CK[z]$.
\end{abstract}

\maketitle

\newcounter{bean}

Given a rational function $\phi\in K(z)$ with coefficients
in a field $K$,
one may consider the dynamical system given by the action
of the iterates $\phi^n$ on $\PK = K\cup\{\infty\}$,
for $n\geq 0$.  Here,
$\phi^n$ denotes the $n$-fold composition $\phi\circ\cdots\circ\phi$,
so that $\phi^0$ is the identity function, $\phi^1=\phi$,
$\phi^2=\phi\circ\phi$, and so on.  The case of complex
dynamics, when $K=\CC$, has been studied intensively for
several decades; see \cite{Bea,CG,Mil} for expositions.
In particular, one may study the action of $\{\phi^n\}$
on the Fatou set $\Fat=\Fat_{\phi}$, to be defined below.  It is
well known that the connected components of $\Fat$
are mapped onto one another by $\phi$ and that, according
to Sullivan's deep No Wandering Domains Theorem \cite{Sul},
every complex Fatou component is preperiodic under application
of $\phi$.

It is also possible to define Fatou sets for other metric
fields. The study of the resulting dynamics has seen growing
interest in the past decade or two;
see, for example,
\cite{Ben2,Bez,HerYoc,Hs2,Lub,Riv1}.
In this paper we will study wandering domains over certain
non-archimedean fields, and we fix the following notation.
\begin{tabbing}
$\CK$ \= \hspace{1.0in} \= a complete and algebraically closed
  non-archimedean field \\
$|\cdot|$ \> \> the absolute value on $\CK$ \\
$\hat{k}$ \> \> the residue field of $\CK$ \\
$p$ \> \> the residue characteristic $\charact\hat{k}$ \\
$\PCK$ \> \> the projective line $\CK\cup\{\infty\}$
\end{tabbing}
We will assume throughout that $p>0$; that is, that
$\CK$ has {\em positive residue characteristic}.
Note that $0\leq |p|<1$, when $p$ is viewed as
an element of $\CK$ under the unique nontrivial homomorphism
of $\ZZ$ into $\CK$.
It is possible that
$\charact\CK=0$ with $\charact\khat=p>0$;
but if $\charact\CK>0$, then $\charact\khat=\charact\CK$.

Recall that ``non-archimedean'' means
$\CK$ satisfies the ultrametric triangle inequality
$$|x+y|\leq\max\{ |x|,|y| \}
\qquad \text{for all} \qquad x,y\in\CK .
$$
If $|x|\neq |y|$, it is immediate that
$|x+y|=\max\{ |x|,|y| \}$.
Note that $|n|\leq 1$ for all $n\in\ZZ$.
Recall also that the residue field $\khat$ is defined to
be $\ints_K/\maxid_K$, where
$\ints_K$ is the ring $\{x\in \CK : |x|\leq 1\}$ of integers
in $\CK$, and
$\maxid_K$ is the maximal ideal of $\{x\in K : |x|< 1\}$ of $\ints_K$.
It is easy to check that
$\khat$ is algebraically closed because $\CK$ is.

The best known example of such a field is $\Cp$, constructed
as follows.  Given a prime integer $p\geq 2$,
the (complete non-archimedean) field $\Qp$ of $p$-adic rational numbers
is the completion of the usual rational numbers $\QQ$ with
respect to the $p$-adic absolute value given by $|p^e r|=p^{-e}$,
where $r$ is a rational number with numerator and denominator
both prime to $p$.
Its algebraic closure is $\Qpbar$, and the completion of
$\Qpbar$ is $\Cp$; the absolute value extends uniquely to
$\Qpbar$ and hence to $\Cp$.
The residue field is $\hat{k}=\Fpbar$, the algebraic closure
of the field $\FF_p$ of $p$ elements.
Note that $\charact \Cp = \charact \Qp=0$, while
$\charact \hat{k} = p>0$.

As another example,
if $L$ is any abstract field, then $K=L((T))$, the field of formal
Laurent series with coefficients in $L$, is a complete non-archimedean
field.
(We may define an absolute value $|\cdot|$ on $K$
by $|f| = 2^{-n}$, where $n\in\ZZ$ is the smallest
integer for which the $T^n$ term of the formal Laurent
series $f$ has a nonzero coefficient.)
Once again, the absolute value extends uniquely to an
algebraic closure of $K$.  If we denote by $\CK$
the completion of an algebraic closure of $K$,
then $\CK$ is algebraically closed and complete;
its residue field $\khat$ is an algebraic closure of $L$.
In this case, $\charact \CK=\charact \khat = \charact K=\charact L$.
We refer the reader to \cite{Kob,Rob,Ser1}
for surveys of non-archimedean fields.

In complex dynamics, given a rational function $\phi\in\CC(z)$
with complex coefficients, one considers
the action of the iterates $\{\phi^n\}$ on the Riemann
sphere $\PC$.  The {\em Fatou set} $\Fat$ of $\phi$ is defined
to be the the set of all points $x\in\PC$
such that the family $\{\phi^n:n\geq 0\}$ is equicontinuous at $x$,
with respect to the spherical metric on $\PC$;
the {\em Julia set} $\Jul$ is the complement
$\PC\setminus\Fat$.
(Recall that a family $F$ of functions from
a metric space $X$ to a metric space $Y$
is called {\em equicontinuous} at $x_0\in X$ if
for every $\eps>0$ there is a $\delta>0$ such that
$d_Y(f(x),f(x_0))<\eps$ for all $f\in F$ and for all $x\in X$
satisfying $d_X(x,x_0)<\delta$.  The key point is that $\eps$ is
chosen independent of $f$.)
Intuitively, the Fatou set is the region of order, where
the iterates $\{\phi^n\}$ are well behaved; the Julia
set is the region of chaos, where a small error becomes
huge after many iterations of $\phi$.

Alternately, by the Arzel\`{a}-Ascoli Theorem, one may define
the Fatou set of a complex function
to be the set of points where $\{\phi^n\}$
forms a normal family.
However, because $\CK$ is not locally compact, and therefore
the Arzel\`{a}-Ascoli Theorem fails, non-archimedean
Fatou sets are usually defined in terms of equicontinuity.

It is easy to verify that the Fatou set is open and that $\phi(\Fat)=\Fat$;
similarly, the Julia set is closed, and $\phi(\Jul)=\Jul$.
One may partition the (complex) Fatou set of $\phi\in\CC(z)$
into connected
components; the function  $\phi$ then maps each component into
(and in fact onto) a component.  Thus, we may speak of fixed,
periodic, and preperiodic components.  A component which is
not preperiodic (that is, a component with infinite forward
orbit) is called a {\em wandering domain} of $\phi$.  In 1985,
Sullivan \cite{Sul} proved, using quasiconformal conjugations,
that a complex rational function $\phi\in\CC(z)$ cannot have
wandering domains.

A similar theory of components of the Fatou set
exists for non-archimedean
rational functions $\phi\in\CK(z)$; see
Section~\ref{sect:fat} below, where we give a
precise characterization of such components for our
setting.  However, in spite
of the existence of certain non-archimedean
versions of Teichm\"{u}ller spaces \cite{Moc}, 
there is no apparent analogue of actual
quasiconformal maps suitable for use as conjugating
functions in $\CK(z)$.  In particular, the author's
constructions of wandering domains over
any algebraically closed complete non-archimedean field
(see \cite{Ben7} for residue characteristic zero,
\cite{Ben6} for $\Cp$, and this paper
for general fields of positive residue characteristic)
suggest that Sullivan's proof has no analogue
for $\CK$, and that no appropriate versions
of quasiconformal maps exist in that context.

Nevertheless, the author showed in \cite{Ben3} that
given some weak hypotheses on $\phi\in\Cp(z)$ (including
the assumption that the coefficients of $\phi$ are
algebraic over $\Qp$), the Fatou set has no wandering
domains.  For more general fields, the arguments
in \cite{Ben2} imply that if $\phi\in\CK(z)$ is defined
over a locally compact subfield of $\CK$, and if the Julia
set of $\phi$ contains no critical points, then $\phi$
has no wandering domains.

We suspect that the condition concerning critical points
is unnecessary.  However,
the condition that $\phi$ be defined over a locally
compact field is crucial.  Indeed,
the main result of this paper is the following theorem.
It shows that if the residue field $\khat$ has
positive characteristic, then
there are many polynomials defined over the
(not locally compact) field $\CK$ which have wandering
domains in $\CK[z]$.

\begin{thmblank}
Let $\CK$ be an algebraically closed complete
non-archimedean field with residue field $\khat$
and residue characteristic $\charact\khat =p>0$.
Let $a_0\in\CK$ with $|a_0|>1$,
and let $\eps>0$.  There
is some $a\in\CK$ with $|a-a_0|\leq\eps$ such that the
Fatou set $\Fat$ of the function
$$\phi_a(z) = (1-a)z^{p+1} + a z^p$$
has a wandering component.
\end{thmblank}

A weaker version of the above result appeared in \cite{Ben6},
with some details of the proof omitted.
There, it was assumed that $\CK=\Cp$, and it was only shown
that there is at least one $a\in\Cp$ giving a wandering domain.
Indeed, an examination of that proof reveals that the point $a$
constructed has $|a|=|p|^{-(p-1)}$.  Moreover, the proof
of the main theorem of \cite{Ben6} fails for fields $\CK$
of positive characteristic.

Although the statement and proof of the above Theorem
is similar in structure to its predecessor in \cite{Ben6},
it is much stronger.
First, the Lemmas in Section~\ref{sect:param}
of the current paper are more complicated than their
analogues in \cite{Ben6}, in order to apply to the
more general class of fields $\CK$ (including
fields of characteristic $p$, not just residue characteristic $p$).
Second, a slight modification of the main proof allows
the broader statement
that the set of parameters $a$ giving wandering
domains is dense in the large open subset $\{a\in\CK :|a|>1\}$ of
the parameter space.
That large subset cannot be extended;
there can be no wandering domains for $\phi_a$ with $|a|\leq 1$,
because such $\phi_a$ are conjugate to maps of the form
$\psi(z) = z^{p+1} + c z^p$ with $|c|\leq 1$, which have
good reduction and therefore empty Julia set.

Our proof uses some computations specific to this family,
but the general method should apply to a larger class of
functions.  In particular, the main properties of $\phi_a$
required in the proof are that it maps
a disk $U_0$ (in this case, about $0$)
onto itself with degree divisible by $p$,
and another disk $U_1$ (in this case, about $1$)
onto $U_0$ with degree prime to $p$.
In the case that $\CK=\Cp$ has characteristic zero,
Fern\'{a}ndez \cite{Fer} has extended the argument in \cite{Ben6}
to show that there is
an open set in the higher-dimensional parameter space of 
all degree $(p+1)$ polynomials in $\Cp[z]$ for which
a dense subset of parameters exhibit wandering domains.
(She adds our map $\phi_a$ to an arbitrary polynomial $Q$ with coefficients
of sufficiently small absolute value
and shows, by essentially the same argument as in \cite{Ben6},
that $a$ can be
adjusted by an arbitrarily small amount to guarantee that
$Q+\phi_a$ has a wandering domain.)  It should be straightforward
to reproduce Fern\'{a}ndez' argument for the more general case
of $\CK[z]$ using the methods of this paper, where we do not assume,
as Fern\'andez and \cite{Ben6} do, that one may divide by $p$.
However, the proof of our Theorem
is already so heavy with notation
that we prefer to avoid the further tedium of proving
the higher-dimensional density statement.
Still, we hope that the interested reader will be able
to glean enough from our proof and from that of Fern\'andez
to deduce how the
statement can be extended to the higher-dimensional
parameter space.

Some comments are in order regarding the field of definition
of the parameters $a$ in the Theorem.  As previously
mentioned, the results of \cite{Ben2} imply that if $a$
were defined over a locally compact subfield $K\subseteq\CK$,
and if $\phi_a$ had no critical points in its Julia
set (that is, if $\phi_a$ were {\em hyperbolic}),
then $\phi_a$ would have no wandering domains.
However,
the parameter $a$ chosen in our proof cannot lie in a discretely
valued field.  Indeed, the resulting map $\phi_a$ will, after many
iterations, map all points $z$ in the wandering disk $U_0$ to points
$\phi_a^n(z)$ for which the ratio $\log |z|/\log|\phi_a^n(z)|$ of valuations
is a rational number of arbitrarily large denominator.  Such a situation
could not occur if the map $\phi_a$ were defined over a discretely valued
field; in particular, $a$ cannot lie in any locally
compact subfield of $\CK$.
In other words,
the truth of the Theorem above requires
the freedom to choose $a\in\CK$; a condition
like $a\in\Qpbar$ would be too restrictive.

\section{Non-archimedean disks, mapping lemmas, and dynamics}
\label{sect:fat}

We will denote the closed disk of radius $r>0$ about
a point $a\in\CK$ by $\Dbar_r(a)$, and the open disk
by $D_r(a)$.  We recall some basic properties
of non-archimedean disks.  Every disk is both open and closed
as a topological set.
Any point in a disk $U$
is a center, but the radius of $U$ is a well-defined real number
$\rad (U)$, which is the same as the diameter of $U$.  If two
disks in $\CK$ intersect, then one contains the other.

The following lemma describes the action of polynomials on disks.
It is easy to prove using
basic non-archimedean analysis results,
such as Hensel's Lemma and the Weierstrass Preparation Theorem.
(See Lemmas~2.2 and~2.6
of \cite[Section 2]{Ben6a}; see also \cite[Chapter 6]{Rob}
for a broad introduction to such techniques over $p$-adic fields.)

\begin{lemma}
\label{lem:diskim}
Let $U\subseteq\CK$ be a disk with $\rad(U)=r$, and let $f\in\CK[z]$
be a non-constant polynomial.  Then:
\begin{list}{\rm \alph{bean}.}{\usecounter{bean}}
\item
$f(U)$ is a disk.
\item
For any $x,y\in U$,
\begin{equation}
\label{eq:mapbdb}
|f(x) - f(y)| \leq \frac{s}{r} \cdot |x-y|,
\end{equation}
where $s=\rad(f(U))$.
\item
$f:U\rightarrow f(U)$ is bijective if and only if
inequality~\eqref{eq:mapbdb} attains equality
for every $x,y\in U$.
\end{list}
\end{lemma}

The projective line $\PCK$ is equipped with a natural
spherical metric
$$\Delta(x,y) = \frac{|x-y|}{(\max\{|x|,1\})(\max\{|y|,1\})}$$
analogous to the spherical metric on $\PC$.
It is easy to verify that $\Delta(x,y)=|x-y|$ if $|x|,|y|\leq 1$;
that $\Delta(x,y)=|x^{-1}-y^{-1}|$ if $|x|,|y|\geq 1$;
and that $\Delta(x,y)=1$ otherwise.

Let $\phi\in\CK(z)$ be a rational function, and consider the
dynamical system of the iterates $\{\phi^n\}$ acting
on $\PCK$.
As in complex dynamics, we define
the {\em Fatou set} $\Fat$ of $\phi$ to the the set of all points $x\in\PCK$
such that $\{\phi^n:n\geq 0\}$ is equicontinuous at $x$
with respect to $\Delta$;
we define the {\em Julia set} $\Jul$ to be the complement
$\PCK\setminus\Fat$.

Following any of \cite{Ben5,Ben7,Riv1}, it
is possible to partition the Fatou set of $\phi$ into ``components''
using a different definition than the usual topological connected
components.  (The various definitions differ slightly, but they
agree on the question of whether or not wandering domains exist
for a given map.)
Rather than repeating the theoretical definitions,
we refer to Theorem~5.4 of \cite{Ben5},
which provides the following characterization of Fatou components
in our situation.

Specifically, suppose that $\phi$ is a polynomial
of degree at least two and the Julia set of $\phi$ is nonempty;
let $\Fat$ denote the Fatou set of $\phi$.
If $x\in\Fat$ has bounded forward orbit (that is, there exists
$R>0$ such that $|\phi^n(x)|<R$ for all $n\geq 0$), then
the component $V$ of $\Fat$ containing $x$
is the largest closed disk $\Dbar_r(x)$ containing $x$
and contained in $\Fat$.
(At least, that is the case for any of the
definitions appearing in \cite{Ben5}.  For
the definition in \cite{Riv1} and the remaining
definitions in \cite{Ben7}, the component is the
open disk $D_r(x)$ of the same radius.)
Moreover, for such $x$, the component of $\phi^n(x)$ is precisely
$\phi^n(V)$.

We recall Hsia's criterion \cite{Hs2} for
equicontinuity, which is a non-archimedean analogue
of the Montel-Carath\'{e}odory Theorem.
\begin{lemma}
\label{lem:hsia}
{\rm (Hsia)} 
Let $F\subseteq\CK[z]$ be a family of polynomials,
and let $U\subseteq \CK$ be a disk.
Suppose that there is a point $y\in\CK$ such that
for all $f\in F$ and $x\in U$, we have
$f(x)\neq y$.  Then $F$ is an equicontinuous family.
\end{lemma}

Hsia stated his result for arbitrary non-archimedean
meromorphic functions; in that setting, the criterion
is that there are two points in $\PCK$ which are
omitted by every $f\in F$.
In our simpler rephrasing above, the first point is $y$,
and the second is $\infty$.  Lemma~\ref{lem:hsia}
follows easily from Lemma~\ref{lem:diskim}
by first making a change of
coordinates to move $y$ to $0$.  Any $f(U)$ is a disk
not containing $0$ (by Lemma~\ref{lem:diskim}.a and the
hypotheses); therefore it either does not intersect
$\Dbar_1(0)$ or is contained in $\Dbar_1(0)$.  If the
former, then $f(U)$ has radius at most $1$.  If the
latter, then the image of $f(U)$ under $1/z$ is a disk
of radius at most $1$.  Because $z\mapsto 1/z$
preserves the spherical metric,
equicontinuity follows from Lemma~\ref{lem:diskim}.b.

\section{The family}
\label{sect:loclem1}

Recall that $p=\charact\khat >0$.
For any $a\in\CK$, we define
$$\phi_a(z) = (1-a)z^{p+1} + a z^p$$
as in the Theorem.  Note that $0$ and $1$ are
fixed points of $\phi_a$.  The point at $0$
is superattracting (meaning that $\phi_a'(0)=0$,
and therefore nearby points are strongly attracted
to $0$ under iteration), and if $|a|>1$, then
the point at $1$ is repelling (meaning
that $|\phi_a'(1)|>1$, and therefore nearby points
are pushed away from $1$ under iteration).

The following
lemma gives a fairly accurate estimate for
the expansion or contraction of distances
under a single application of $\phi_a$.

\begin{lemma}
\label{lem:distort}
Let $a\in\CK$ with $|a|>1$, and let $y_1,y_2\in\CK$ with
$|y_1|\geq |y_2|$.  Then
$$ |\phi_a(y_1) -\phi_a(y_2)| \leq
	|y_1 - y_2| \cdot |a| \cdot |y_1|^{p-1} \cdot
	\max\left\{ |p|, |y_1|,
	\left( \frac{|y_1 -y_2|}{|y_1|}\right)^{p-1} \right\}.$$
Furthermore, if $y_1,y_2\in D_1(1)$, then
$|\phi_a(y_1) -\phi_a(y_2)| = |a|\cdot |y_1 - y_2|$.
\end{lemma}

{\bf Proof of Lemma~\ref{lem:distort}.}
We have
$$\phi_a(y_1) -\phi_a(y_2) = 
	(1-a) (y_1^{p+1} - y_2^{p+1}) + a (y_1^p - y_2^p) $$
\begin{equation}
\label{eq:expand}
= a (y_1 - y_2) \cdot \left[
	\left(\frac{1}{a}-1\right)
	\left( \sum_{j=0}^{p} y_1^j y_2^{p-j} \right)
	+ (y_1-y_2)^{p-1}
	+ p \left( \sum_{j=0}^{p-1} A_j y_1^j y_2^{p-1-j} \right) \right],
\end{equation}
with $A_j\in\ZZ$.  Thus,
$$|\phi_a(y_1) -\phi_a(y_2)| \leq
	|y_1 - y_2| \cdot |a| \cdot
	\max\left\{ |y_1|^p, |y_1 - y_2|^{p-1}, |p|\cdot |y_1|^{p-1}
	\right\},$$
as claimed.  Finally, if $y_1,y_2\in D_1(1)$, then
$$\sum_{j=0}^{p} y_1^j y_2^{p-j} \in D_1(1) \quad\text{and}\quad
	|y_1 -y_2|<1.$$
In that case, then, \eqref{eq:expand} has absolute value
$|a|\cdot |y_1 - y_2|$.
\qed

\section{Mapping Properties of the Iterates}
\label{sect:loclem2}

In this section, we will
apply Lemma~\ref{lem:distort} to describe the
behavior of iterates $\phi_a^i$
near the fixed points at $z=0,1$, for $i$ relatively small.

Fix $a_0\in\CK$ with $|a_0| > 1$
(as chosen in the hypotheses of the Theorem),
and to simplify future notation, define
\begin{equation}
\label{eq:mudef}
R = |a_0|^{-1/(p-1)},
\qquad
\mu = \max\{|p|,R\}, \quad \text{and} \qquad
S =  \mu R^3.
\end{equation}
Note that $0<S<R^{-2}S<R\leq \mu <1$.
The quantity $R$ represents the radius
of a particular disk about $0$.  Specifically,
for $a\in\CK$ with $|a|=|a_0|>1$,
$$\phi_a(\Dbar_R(0)) = \Dbar_R(0);$$
the mapping is $p$-to-$1$.
The smaller quantity $S$
is important because $\phi_a^i$ acts
in a predictable manner
on disks of radius at most $S$
just outside $\Dbar_R(0)$, as shown in the following lemma.

\begin{lemma}
\label{lem:shrdisk}
Let $a\in\CK$ with $|a|=|a_0|$, and let $m\geq 1$.
Let $x\in\CK$ with $\dsps R<|x|\leq R^{1-p^{-m}}$.
Then for all $0\leq i\leq m$, we have
$$|\phi_a^i(x)| = R^{1-p^i}|x|^{p^i}\leq 1,
\qquad\text{and}\qquad
\phi_a^i(\Dbar_S(x))\subseteq \Dbar_{S'_i}(\phi_a^i(x)),$$
where
$$S'_i = \mu^i R^{-e_i} S < \mu^i R^{-2} S,
\qquad\text{and}\qquad
e_i = p^{1-m} + p^{2-m} + \cdots + p^{i-m}.$$
\end{lemma}

{\bf Proof of Lemma~\ref{lem:shrdisk}.}
Observe that $e_i$ is a partial sum of a geometric series,
so that $e_i < p/(p-1) \leq 2$.  In particular, the inequality
$S'_i < \mu^i R^{-2} S$ is valid.
Meanwhile, for any $y\in D_1(0)$, 
$|\phi_a(y)|=|a|\cdot |y|^p$,
and so the first statement follows from a simple induction.

We now prove the second statement of the lemma by induction on $i$.
If $i=0$, in which case $e_i=0$, the statement is clear.
Given that the statement holds for some $0\leq i\leq m-1$,
write $y=\phi_a^i(x)$ to simplify notation.
Note that $S'_i=\mu^{i+1} R^{3-e_i} < \mu R$,
and that
$$\frac{(\mu R)^{p-1}}{|y|^p} <
\frac{\mu}{|y|} \leq
\max\left\{\frac{|p|}{|y|} , 1\right\}.$$
Therefore, by Lemma~\ref{lem:distort},
\begin{equation}
\label{eq:subset}
\phi_a\left(\Dbar_{S'_{i}}(y)\right)
\subseteq \Dbar_{\lambda(y) \cdot S'_{i}}(\phi_a(y)),
\qquad\text{where}\qquad
\lambda(y) = |a|\cdot |y|^p \max\left\{1,\frac{|p|}{|y|}\right\}.
\end{equation}
We compute
\begin{align*}
\lambda(y) & = \lambda\left(\phi_a^i(x)\right) =
  |a| \cdot \left| \phi_a^i(x) \right|^p \cdot
  \max\left\{ 1, \left| \frac{p}{\phi_a^i(x)}\right| \right\}
  \\
  & \leq |a| \left( R^{1-p^i} |x|^{p^i} \right)^p
  \cdot \max\left\{ 1, \frac{|p|}{R} \right\}
  \\
  & = |a| R^{p-1-p^{i+1}} |x|^{p^{i+1}} \mu
  = \mu R^{-p^{i+1}} |x|^{p^{i+1}}
  \leq \mu R^{-p^{i+1-m}} .
\end{align*}
Thus, $\lambda(y)\cdot S'_i \leq S'_{i+1}$.
By the inductive hypothesis,
$\phi_a^i(\Dbar_S(x))\subseteq \Dbar_{S'_i}(y)$;
hence, by \eqref{eq:subset},
$$\phi_a^{i+1}(\Dbar_S(x))\subseteq \Dbar_{S'_{i+1}}(\phi_a(y)),$$
and the Lemma follows.
\qed

\begin{remark}
{\em a.} 
The upper bound of $S'_i$ for the image radius in Lemma~\ref{lem:shrdisk}
is far from sharp.  The sharp bound is not difficult to compute, but
our Theorem will require only a bound
which decreases to $0$ as $i$ increases; the simpler bound
given above will suffice.

\begin{figure}[b]
\scalebox{1.13}{\includegraphics{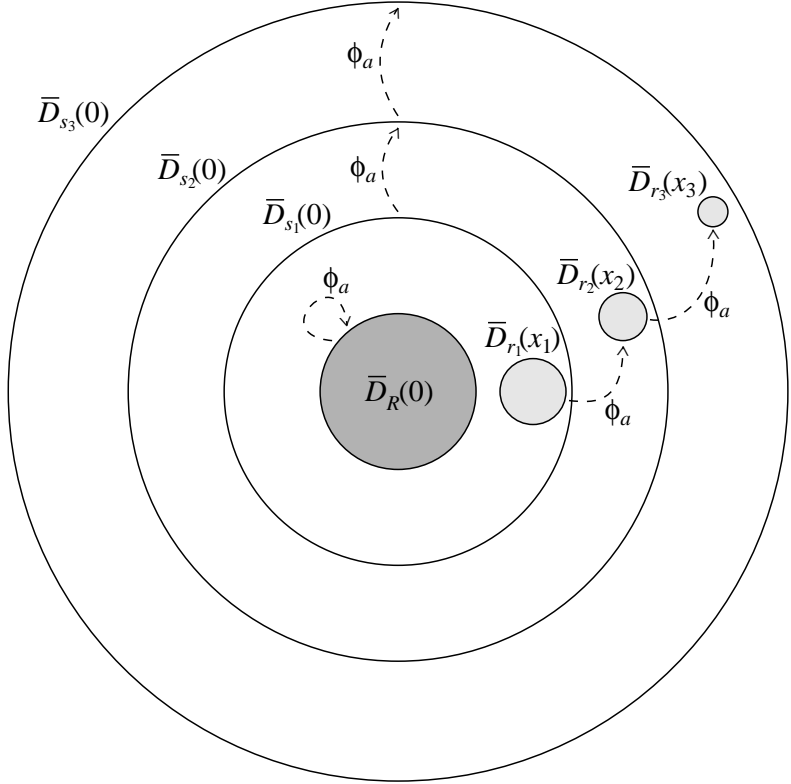}}
\caption{Global expansion and local contraction under $\phi_a$.}
\label{fig:wanddom}
\end{figure}

{\em b.}
Intuitively, Lemma~\ref{lem:shrdisk} says the following.  A point $x$
just outside of $\Dbar_R(0)$ is pushed further away from $0$ by $\phi_a$.
However, $|\phi_a'(x)| < 1$, so that locally
(i.e., on $\Dbar_S(x)$), $\phi_a$ is contracting.
The situation is illustrated in Figure~\ref{fig:wanddom},
where the small disk $\Dbar_{r_1}(x_1)$
maps to the smaller disk $\Dbar_{r_2}(x_2)$ and then to the even
smaller disk $\Dbar_{r_3}(x_3)$.  However, at the same time, with
$s_i=|x_i|$, the disk $\Dbar_{s_1}(0)$ maps to the larger disk
$\Dbar_{s_2}(0)$ and then to the even larger disk $\Dbar_{s_3}(0)$.
In broad terms, $\phi_a$ is locally contracting on the whole
annulus $R<|x|<1$, even as every point in the annulus is
pushed away from $0$.  Thus, if $|x|$
is only slightly bigger than $R$, the iterates of $x$ will stay
inside the larger disk $D_1(0)$ for some time; and the corresponding
iterates of $\Dbar_S(x)$ will become very small in radius.

{\em c.}
Ultimately, our main interest in Lemma~\ref{lem:shrdisk}
lies in the case that $i=m$.
However, the cases $i<m$ are important both for the
inductive proof and because we need to have some idea of the
orbit traversed by $x$ during those $m$ iterations.
An analogous remark applies to the quantity $M$ in
Lemma~\ref{lem:repdisk} below.
\end{remark}

\begin{lemma}
\label{lem:repdisk}
Let $a\in\CK$ with $|a|=|a_0|$, and let $M\geq 1$.
Let $x\in\CK$ with $|x-1|\leq |a|^{-M}$.
Then for all $0\leq i\leq M$, we have
$$|\phi_a^i(x)-1| = |a|^i\cdot |x-1|,$$
and for any $0<r\leq |a|^{-M}$,
$$\phi_a^i(\Dbar_r(x)) = \Dbar_{|a|^i r}(\phi_a^i(x)).$$
\end{lemma}

{\bf Proof of Lemma~\ref{lem:repdisk}.}
The result follows easily by induction on $i$ using
the first statement of Lemma~\ref{lem:distort}
and the fact that $\phi_a(1)=1$.
\qed

\section{Varying the parameter}
\label{sect:param}

We wish to study the effects of perturbations of $a$
on the dynamics of $\phi_a$, and so we introduce
the following notation.
For any $z\in\PCK$, $a\in\CK$, and $n\geq 0$, define
\begin{equation}
\label{eq:Phidef}
\Phi_n(a,z) = \phi_a^n(z).
\end{equation}
For any fixed $x\in\CK$,
$\Phi_n(a,x)$ is a polynomial in the variable $a$.
In particular, for any disk $D\in\CK$, the set
$\Phi_n(D,x)$ is also a disk, by Lemma~\ref{lem:diskim}.a.

\begin{lemma}
\label{lem:repwvar}
Let $a\in\CK$ with $|a|=|a_0|$, let $M\geq 0$, let $n\geq 1$, and let
$x\in\CK$ with $|\phi_a^n(x)-1|\leq |a|^{-M}$.
Let $\eps\in (0,|a|^{1-M}]$.
Suppose that $\Phi_n(\cdot, x)$ maps
$\Dbar_{\eps}(a)$ bijectively onto
$\Dbar_{\eps/|a|}(\phi_a^{n}(x))$.  Then
$\Phi_{n+M}(\cdot,x)$ maps $\Dbar_{\eps}(a)$ bijectively
onto $\Dbar_{\eps\cdot |a|^{(M-1)}}(\phi_a^{n+M}(x))$.
\end{lemma}

{\bf Proof.}
We proceed by induction on $M$; the case $M=0$ is vacuously true.
For $M\geq 1$,
let $b_1,b_2\in\Dbar_{\eps}(a)$,
let $y_1=\Phi_{n+M-1}(b_1,x)$, and let $y_2=\Phi_{n+M-1}(b_2,x)$.
Note that
\begin{align*}
|y_1 - 1| & = |a|^{M-1} \cdot \left| \Phi_n(b_1,x) - 1 \right| \\
	& \leq |a|^{M-1} \cdot
	\max\left\{ \left| \Phi_n(b_1,x) - \phi_a^n(x)\right| ,
	\left| \phi_a^n(x) - 1 \right| \right\} \\
	& \leq |a|^{M-1} \cdot
	\max\left\{ |a|^{-1}\cdot |b_1-a|, |a|^{-M} \right\}
	= |a|^{-1},
\end{align*}
where the first equality is by Lemma~\ref{lem:repdisk},
and the rest is by the hypotheses, ultrametricity,
and Lemma~\ref{lem:diskim}.b.
Therefore,
\begin{align*}
\left| \phi_{b_1}(y_1) - \phi_{b_2}(y_1) \right| & =
  \left| (1-b_1) y_1^{p+1} + b_1 y_1^p - 
  (1-b_2) y_1^{p+1} - b_2 y_1^p  \right|
  \\
  & = |y_1|^p \cdot |y_1 - 1| \cdot |b_1-b_2| 
  \leq |a|^{-1}\cdot |b_1 - b_2|,
\end{align*}
by the above bound and because $|y_1|=1$.
On the other hand, by Lemma~\ref{lem:distort},
$$\left| \phi_{b_2}(y_1) - \phi_{b_2}(y_2) \right| =
	|a|\cdot |y_1 - y_2| = |a|^{M-1}\cdot |b_1-b_2|,$$
where the final equality is by induction
and Lemma~\ref{lem:diskim}.c.
In particular, if $b_1\neq b_2$,
$$\left| \phi_{b_1}(y_1) - \phi_{b_2}(y_1) \right| <
\left| \phi_{b_2}(y_1) - \phi_{b_2}(y_2) \right|,$$
and therefore
$$\left| \phi_{b_1}(y_1) - \phi_{b_2}(y_2) \right| =
\left| \phi_{b_2}(y_1) - \phi_{b_2}(y_2) \right|$$
by ultrametricity.  It follows that
\begin{equation}
\label{eq:nM}
\left| \Phi_{n+M}(b_1,x) - \Phi_{n+M}(b_2,x) \right|
	= |a|^{M-1}\cdot |b_1-b_2|.
\end{equation}
The same equality is clear if $b_1=b_2$; thus,
\eqref{eq:nM} holds for all $b_1,b_2\in\Dbar_{\eps}(a)$.

By Lemma~\ref{lem:diskim}.a,
$\Phi_{n+M}(\Dbar_{\eps}(a),x)$ is a disk about $\phi_a^{n+M}(x)$
of some radius $s$.  By \eqref{eq:nM} and Lemma~\ref{lem:diskim}.c,
$s=\eps \cdot |a|^{(M-1)}$, and the mapping is bijective, as claimed.
\qed

\begin{lemma}
\label{lem:shrwvar}
Let $a\in\CK$ with $|a|=|a_0|$, let $m\geq 1$, let $n\geq 0$, and let
$x\in\CK$ with $|\phi_a^n(x)|=R^{1-p^{-m}}$.  
Let $A>|a|^{-1}$ be a real number, let $\eps\in (0,A^{-1}S]$,
and suppose that
$$\Phi_n(\Dbar_{\eps}(a),x)\subseteq
	\Dbar_{A\eps}(\phi_a^n(x))
	\quad\text{and}\quad
	\mu^m < A^{-1}R^{p+1}.$$
Then $\Phi_{n+m}(\cdot,x)$ maps $\Dbar_{\eps}(a)$ bijectively
onto $\Dbar_{\eps/|a|}(\phi_a^{n+m}(x))$.
\end{lemma}

{\bf Proof.}
For any $0\leq i \leq m$, let
$e_i = p^{1-m} + p^{2-m} + \cdots + p^{i-m}$
as in Lemma~\ref{lem:shrdisk},
and let
$$r_i = \max\left\{ \mu^i R^{-e_i} A, R^{p-p^{-(m-i)}}\right\} .$$
Note that $r_i \leq A$, and therefore $r_i \eps \leq A\eps \leq S$.
We claim that for any such $i$,
\begin{equation}
\label{eq:ribound}
\Phi_{n+i}\left( \Dbar_{\eps}(a),x \right) \subseteq
\Dbar_{r_i \eps}\left( \phi^{n+i}_a(x) \right),
\end{equation}
with $\Phi_{n+i}(\cdot, x)$ mapping $\Dbar_{\eps}(a)$
bijectively onto $\Dbar_{r_i\eps }( \phi^{n+i}_a(x) )$ if
\begin{equation}
\label{eq:mucond}
i\geq 1 \quad \text{and} \quad
\mu^i R^{-e_i} A < R^{p-p^{-(m-i)}}.
\end{equation}
Condition~\eqref{eq:mucond} holds for $i=m$ because
$\mu^m < A^{-1}R^{p+1}$ by hypothesis and because $e_m<2$;
moreover, $r_m= R^{p-1}= |a|^{-1}$.
Thus, to prove the lemma, it suffices to prove the claim.

We proceed by induction on $i$.  Since
$r_0=A$,  the $i=0$ case is true by hypothesis.
For $1\leq i \leq m$, pick any two distinct points
$b_1,b_2\in \Dbar_{\eps}(a)$,
and let $y_1=\Phi_{n+i-1}(b_1,x)$ and $y_2=\Phi_{n+i-1}(b_2,x)$.
By Lemma~\ref{lem:diskim}.b, the inductive hypothesis implies
that $|y_1-y_2|\leq r_{i-1}\cdot |b_1 - b_2|$.
Note that $\dsps |y_1|=|y_2|=R^{1-p^{-(m-i+1)}}$, by
Lemma~\ref{lem:shrdisk}.  Therefore,
\begin{align*}
\left| \phi_{b_1}(y_1) - \phi_{b_2}(y_1) \right| & =
  \left| (1-b_1) y_1^{p+1} + b_1 y_1^p - 
  (1-b_2) y_1^{p+1} - b_2 y_1^p  \right|
  \\
  & = |b_1 - b_2| \cdot |y_1 - 1| \cdot |y_1|^p
  = |b_1 - b_2|\cdot |y_1|^p
  \\
  & = |b_1 - b_2|\cdot R^{p-p^{-(m-i)}}
  \leq r_i |b_1 - b_2|,
\end{align*}
with equality if condition~\eqref{eq:mucond} holds.
Meanwhile, $|y_1-y_2|\leq r_{i-1}\eps \leq S<\mu R$;
by Lemma~\ref{lem:distort},
$$\left| \phi_{b_2}(y_1) - \phi_{b_2}(y_2) \right| \leq
	\lambda(y_1) \cdot |y_1 - y_2| \leq
	r_{i-1} \cdot \lambda(y_1) \cdot |b_1 - b_2| ,$$
where $\lambda(y_1) = |a| \cdot |y_1|^p \max\{1, |p|/|y_1|\}$
as in the proof of Lemma~\ref{lem:shrdisk}.
Adding, we have
$$\left| \Phi_{n+i}(b_1,x) - \Phi_{n+i}(b_2,x) \right| =
        \left| \phi_{b_1}(y_1) - \phi_{b_2}(y_2) \right| \leq
        \max\{r_i , r_{i-1} \cdot \lambda(y_1) \} \cdot |b_1 - b_2| ,$$
with equality if the maximum is uniquely attained.
Thus, to prove~\eqref{eq:ribound}, it suffices to show
$$\lambda(y_1) \cdot r_{i-1} \leq r_i.$$
By Lemma~\ref{lem:diskim},
the full claim would follow if in fact
$\lambda(y_1) \cdot r_{i-1} < r_i$
whenever condition~\eqref{eq:mucond} holds.

Note that
$$\lambda(y_1)= R^{-p^{-(m-i)}}\cdot
	\max\left\{ R, |p|\cdot R^{p^{-(m-i+1)}} \right\}
	\leq \mu R^{-p^{-(m-i)}}.$$
Thus, if $r_{i-1}=\mu^{i-1} R^{-e_{i-1}} A$, then
$$\lambda(y_1) \cdot r_{i-1} \leq \mu^i R^{-e_i} A\leq r_i,$$
with strict inequality if condition~\eqref{eq:mucond} holds.
Similarly, if $r_{i-1}= R^{p-p^{-(m-i+1)}}$, then
$$\lambda(y_1) \cdot r_{i-1} =
	R^{p-p^{-(m-i)}}\cdot
	\max\left\{ R^{1-p^{-(m-i+1)}}, |p|\right\}
	< R^{p-p^{-(m-i)}}\leq r_i.$$
The claim is proved, and the lemma follows.
\qed

\begin{remark}
{\em a.}
The two conditions displayed in the hypotheses of
Lemma~\ref{lem:shrwvar} say, first, that $A$ is large enough
to measure the size of a certain image disk, and second, that
$m$ is large enough that $\mu^m$ is smaller than the
already small $A^{-1}R^{p+1}$.

{\em b.}
Lemmas~\ref{lem:repwvar} and~\ref{lem:shrwvar} both
describe the behavior of iterates $\phi_a^{n+i}(x)$ as $a$ varies;
the former is for $\phi_a^n(x)$ close to $1$, and the latter is
for $\phi_a^n(x)$ just outside $\Dbar_R(0)$.
The conclusion of Lemma~\ref{lem:shrwvar} exactly
matches the key hypothesis of Lemma~\ref{lem:repwvar},
and the conclusion of Lemma~\ref{lem:repwvar} implies
the key hypothesis of Lemma~\ref{lem:shrwvar}.  (We
shall need the full strength of Lemma~\ref{lem:repwvar}
elsewhere in the main proof.)
Thus, we will be able to alternate applications
of each Lemma in order to describe the behavior of $\phi_a^n$
for larger and larger values of $n$.
\end{remark}

\section{Proof of the Theorem}

To prove the Theorem, we will first choose several auxiliary
sequences $\{m_i\}$, $\{M_i\}$, $\{n_i\}$ and $\{N_i\}$
of positive integers, and $\{\eps_i\}$ and $\{r_i\}$
of positive real numbers.
The carefully chosen integer sequences will be used to describe
the orbit of a certain point $x$ under various
maps $\phi_a$; the real values $\eps_i$
will be radii of disks in which the parameters $a$ lie;
and the $r_i$ will be radii of disks in which certain
images $\phi_a^n(x)$ lie.

After those choices, we will argue inductively
using Lemmas~\ref{lem:repwvar} and~\ref{lem:shrwvar}
to produce a sequence $\{a_i\}$ of parameters with
limit $a$ such that $\phi_a$ has a wandering domain.

\subsection{Auxiliary sequences}
We are given $\eps >0$ and $a_0\in\CK$ with $|a_0|>1$.
Define $R$, $S$, and $\mu$ as in \eqref{eq:mudef}.
We assume without loss that $\eps\leq S$.

Choose a strictly
increasing sequence of positive integers $\{M_i\}_{i\geq 1}$ so that
$$ |a_0|^{1-M_1}\leq\eps , \qquad
  \text{and for all } i\geq 2, \quad
  |a_0|^{M_{i-1} - M_i} \leq S.$$
Let $\eps_0 = \eps$, and for all $i\geq 1$, let
$$\eps_i = |a_0|^{1-M_i} S.$$
Then $\{\eps_i\}$ is a strictly decreasing sequence
of positive real numbers with limit $0$, and
$|a_0|^{1-M_i} \leq \eps_{i-1}$ for all $i\geq 1$.

Choose another strictly increasing sequence of positive
integers $\{m_i\}_{i\geq 0}$ so that
$$\mu^{m_i}\leq |a_0|^{-M_{i+1}} \cdot R^2
	\quad\text{for all } i\geq 0.$$
For all $i\geq 0$, define
$$r_i=R^{1-p^{-m_i}}.$$
Then $\{r_i\}$ is a strictly decreasing sequence of
positive real numbers with limit $R$.
Note that the starting index for $\{m_i\}$, $\{\eps_i\}$,
and $\{r_i\}$ is $i=0$, while
that for $\{M_i\}$ is $i=1$.

We will find $x\in D_1(0)$ and $a\in\Dbar_{\eps}(a_0)$
so that the orbit
$\{\phi_a^j(x)\}_{j\geq 0}$
is described by
\begin{equation}
\label{eq:orbit}
\underbrace{0,\ldots,0}_{m_0},
\underbrace{1,\ldots,1}_{M_1},
\underbrace{0,\ldots,0}_{m_1},
\underbrace{1,\ldots,1}_{M_2},
\underbrace{0,\ldots,0}_{m_2},
\ldots
\end{equation}
where a $0$ in the $j$-th position in the sequence indicates that
$\phi_a^{j-1}(x)\in D_1(0)$, and a $1$ indicates that
$\phi_a^{j-1}(x)\in D_1(1)$.

To simplify future notation,
define
$$n_i = \sum_{j=1}^{i} (m_{j-1} + M_j) ,
	\quad \text{ and } \quad
N_i = n_i + m_i$$
for all $i\geq 0$.
That is, $n_i$ is the number of terms
in \eqref{eq:orbit}
up to but not including the block of $m_i$ $0$'s, and
$N_i$ is the number of terms
up to but not including the block of $M_{i+1}$ $1$'s.

\subsection{Choosing the parameter}
Recall that $r_0=R^{1-p^{-m_0}}$.
Pick $y\in\CK$ with $|y|=r_0$;
by Lemma~\ref{lem:shrdisk},
$|\phi_{a_0}^{m_0}(y)| = 1$.
(Such a $y$ exists because $r_0\in\CK$.)
In addition, $\phi_{a_0}^{m_0}(0)=0$.
Therefore, by Lemma~\ref{lem:diskim}.a,
$\phi_{a_0}^{m_0}\left(\Dbar_{r_0}(0)\right)$
is a disk containing both $0$ and $\phi_{a_0}^{m_0}(y)$,
and hence
$$\Dbar_1(0)\subseteq \phi_{a_0}^{m_0}\left(\Dbar_{r_0}(0)\right).$$
Thus, there is some
$x\in\Dbar_{r_0}(0)$ with $\phi_{a_0}^{m_0}(x)=1$.
Again by Lemma~\ref{lem:shrdisk}, we have $|x|=r_0$.

For every $i\geq 0$, we will find $a_i\in\Dbar_{\eps_{i-1}}(a_{i-1})$
(or just the original $a_0$, for $i=0$)
so that for every $a\in\Dbar_{\eps_i}(a_i)$,
the orbit $\{\phi_{a}^j(x)\}$ follows
\eqref{eq:orbit} up to the $j=N_i$ iterate,
with $\phi_{a_i}^{N_i}(x)=1$ and such that
\begin{equation}
\label{eq:indcond}
\Phi_{N_i}(\cdot, x) : \Dbar_{\eps_i}(a_i) \rightarrow
	\Dbar_{\eps_i/|a_i|}(1)
	\quad \text{is bijective.}
\end{equation}
Because $\{\eps_i\}$ is decreasing,
each $a_i$ will lie in $\Dbar_{\eps_0}(a_0)$, and
therefore $|a_i|=|a_0|$.

We construct the sequence $\{a_i\}$ inductively.
For $i=0$, we already have
$\phi_{a_0}^{N_0}(x)=1$, and by Lemma~\ref{lem:shrdisk},
the orbit $\{\phi_{a_0}^j(x)\}$ follows
\eqref{eq:orbit} up to the $N_0=m_0$ iterate.  By Lemma~\ref{lem:shrwvar}
(with $n=n_0=0$, $m=m_0$, $a=a_0$, $A=1$, and $\eps=\eps_0$),
condition~\eqref{eq:indcond} holds.  Also, by Lemma~\ref{lem:shrdisk},
the orbit $\{\phi_a^j(x)\}$ is correct up to $j=N_0$
for every $a\in\Dbar_{\eps_0}(a_0)$.
Hence, the $i=0$ case is already done.

For $i\geq 1$, given $a_{i-1}$ with the desired
properties, set $\rho = |a_0|^{1-M_i}\leq \eps_{i-1}$.
By the inductive hyphothesis,
for each $a\in\Dbar_{\rho}(a_{i-1})$,
the orbit $\{\phi_a^j(x)\}$ agrees with~\eqref{eq:orbit}
up to $j=N_{i-1}$.
Moreover, because $\Phi_{N_{i-1}}(\cdot, x)$ maps
$\Dbar_{\eps_{i-1}}(a_{i-1})$ bijectively onto
$\Dbar_{\eps_{i-1}/|a_{i-1}|}(1)$
with $\Phi_{N_{i-1}}(a_{i-1},x)=1$,
we see that $\Phi_{N_{i-1}}(\cdot, x)$ also maps
$\Dbar_{\rho}(a_{i-1})$ bijectively onto
$\Dbar_{\rho/|a_{i-1}|}(1)$,
by Lemma~\ref{lem:diskim}.c.
By Lemma~\ref{lem:repwvar}
(with $a=a_{i-1}$, $M=M_i$, $n=N_{i-1}$, and $\eps=\rho$),
$\Phi_{n_i}(\cdot,x)$ maps $\Dbar_{\rho}(a_{i-1})$
bijectively onto $\Dbar_1(1)$.
Hence,
there exists $b\in\CK$ such that $\Phi_{n_i}(b,x)=0$.
By Lemma~\ref{lem:diskim}.c, we must have $|b-a_{i-1}|=\rho$;
and because $r_i=R^{1-p^{-m_i}}<1$,
\begin{equation}
\label{eq:nimap}
	\Phi_{n_i}\left(\cdot,x \right) :
	\Dbar_{\sigma}(b) \rightarrow \Dbar_{r_i}(0)
	\quad\text{is bijective,}
\end{equation}
where $\sigma = r_i\cdot |a_0|^{1-M_i} = r_i\cdot\rho \in (0,\rho)$.
Moreover, for all $a\in\Dbar_{\sigma}(b)$,
$$|\Phi_{N_{i-1}}(a,x) - 1| = |a_0|^{-1} |a-a_{i-1}|
= |a_0|^{-1} |b-a_{i-1}| = |a_0|^{-M_i}.$$
By Lemma~\ref{lem:repdisk}, then, the orbit
$\{\phi_a^j(x)\}$ is correct up to $j=n_i$ for every
$a\in\Dbar_{\sigma}(b)$.

Choose $c\in\Dbar_{\sigma}(b)$ so that
$|\Phi_{n_i}(c,x)| = r_i$.
(Such $c$ exists because of \eqref{eq:nimap} and the
fact that $r_i \in | \CK |$.)
By Lemma~\ref{lem:shrdisk},
$|\Phi_{N_i}(c,x)|=1$.  Furthermore, it is clear
that $\Phi_{N_i}(b,x)=0$.  Because the polynomial image
of a disk is a disk (Lemma~\ref{lem:diskim}.a), it follows that
$\Phi_{N_i}(\Dbar_{\sigma}(b,x)) \supseteq \Dbar_1(0)$.
We may therefore choose
$a_i\in\Dbar_{\sigma}(b)$ so that $\Phi_{N_i}(a_i,x) = 1$.

By Lemma~\ref{lem:shrdisk}, we must have
$|\Phi_{n_i}(a_i,x)|=r_i$.
Observe also that $\eps_i\leq\sigma$.
By \eqref{eq:nimap} and Lemma~\ref{lem:diskim}.c,
$\Phi_{n_i}(\Dbar_{\eps_i}(a_i),x)$ must be a disk
of radius $\eps_i\cdot |a_0|^{M_i - 1}=S$.
Therefore, by Lemma~\ref{lem:shrwvar}
(with $n=n_i$, $m=m_i$, $a=a_i$,
$A=|a_i|^{M_i - 1}$, and $\eps=\eps_i$),
condition~\eqref{eq:indcond} holds on $\Dbar_{\eps_i}(a_i)$.
By Lemma~\ref{lem:shrdisk},
the orbit $\{\phi_{a}^j(x)\}$ is correct up to $j=N_i$ for every
$a\in\Dbar_{\eps_i}(a_i)$.
Our construction of $a_i$ is complete.

The sequence $\{a_i\}_{i\geq 0}$ is a Cauchy sequence, because
for any $0\leq i \leq j$, we have
$|a_i - a_j|\leq \eps_i,$
and $\eps_i\rightarrow 0$.  Therefore, the sequence has a limit
$a\in\CK$, with $|a-a_0|\leq\eps_0$.
We claim that $\phi_a$ has a wandering domain.

\subsection{The end of the proof}
By construction, $a\in\Dbar_{\eps_i}(a_i)$ for every $i\geq 0$;
hence, the orbit
$\{\phi_a^j(x)\}$ follows \eqref{eq:orbit} exactly.
In light of Lemmas~\ref{lem:shrdisk} and~\ref{lem:repdisk},
we must have
$$|\phi_a^{n_i}(x)|=r_i=R^{1-p^{-m_i}},
\quad\text{and} \quad
|\phi_a^{N_i}(x)-1|=|a|^{-M_{i+1}}$$
for any $i\geq 0$.
Let $U=\Dbar_S(x)$; we will show that $U$ is contained
in a wandering domain of the Fatou set $\Fat$ of $\phi_a$.

Every iterate $U_n=\phi_a^n(U)$ is a disk;
we claim that for any $i\geq 0$,
the radius of $U_{n_i}$ is at most $S$, and the
radius of $U_{N_i}$ is at most $|a|^{-M_{i+1}} S$.
The claim is easily proven by induction, as follows.
For $i=0$, $U_{n_0}=U_0=U$, which has radius $S$;
and by Lemma~\ref{lem:shrdisk},
$U_{N_0}$ has radius at most $\mu^{m_0}R^{-2} S \leq |a|^{-M_1} S$.
For $i\geq 1$, we assume the radius of $U_{N_{i-1}}$
is at most $|a|^{-M_i} S$.  By Lemma~\ref{lem:repdisk},
the radius of $U_{n_i}$ is at most $S$; and
by Lemma~\ref{lem:shrdisk}, the radius
of $U_{N_i}$ is at most $\mu^{m_i}R^{-2} S \leq |a|^{-M_{i+1}} S$.

If some $U_{n_i}$ contained the point $1$, then by
Lemma~\ref{lem:repdisk}, some $U_{n_{i+j}}$ would contain
$\Dbar_1(0)$ and therefore would have radius greater than $S$.
Thus, no $U_{n_i}$ contains the point $1$;
and because
$1$ is fixed, it follows that no $U_n$ contains $1$.
By Hsia's criterion (Lemma~\ref{lem:hsia}), then,
the family $\{\phi_a^n\}$ is equicontinuous on $U$, and
therefore $U\subseteq\Fat$.

The iterates of $U$ are bounded away from $\infty\in\PCK$,
and $\phi_a$ is a polynomial.  Therefore,
by any of the definitions of components in \cite{Ben5,Ben7,Riv1},
the component $V$ of $\Fat$ containing $U$ must be a disk
\cite[Theorem~5.4.d]{Ben5}, as must all of its iterates.
As before, no iterate $\phi^n(V)$
of $V$ can contain the point $1$.  In addition, each $\phi^n(V)$
is contained either in $D_1(0)$ or $D_1(1)$, since it is a disk
intersecting one of those two disks but not containing $1$.
The symbolic dynamics of
$V$ are therefore described by equation~\eqref{eq:orbit}.
Because those
dynamics are not preperiodic, it follows
that $V$ must be wandering.
\qed

\bibliographystyle{plain}

\end{document}